\documentclass[10pt]{article}
\usepackage{amsthm,amsmath,amssymb,amscd}
\usepackage{amssymb}
\usepackage{graphics}

\title{Constant $k$-curvature hypersurfaces in Riemannian manifolds}
\author{Fethi Mahmoudi \thanks{Email: mahmoudi@sissa.it}\\[4mm]
\small  SISSA, Sector of Mathematical Analysis \\ Via Beirut 2-4,
34014 Trieste, Italy}

\date{}

\newtheorem{theorem}{Theorem}[section]
\newtheorem{proposition}{Proposition}[section]

\newcommand{\R}{\mathbb{R}}

\newcommand{\e}{\varepsilon}
\newcommand{\del}{\partial}
\newcommand{\ds}{}

\newcommand{\calO}{{\mathcal O}}

\begin{document}

\maketitle \noindent {\sc abstract}. In \cite{Ye-1}, Rugang Ye
proved the existence of a family of constant mean curvature
hypersurfaces in an $m+1$-dimensional Riemannian manifold
$(M^{m+1},g)$, which concentrate at a point $p_0$ (which is
required to be a nondegenerate critical point of the scalar
curvature), moreover he proved that this family constitute a
foliation of a neighborhood of $p_0$. In this paper we extend this
result to the other curvatures (the $r$-th mean curvature for
$1\le r\le m$). %and we give the expansion of the $m$-dimensional
%volume of the leaves of this foliation as well as the
%$m+1$-dimensional volume of the sets enclosed by each leaf.

\begin{center}

\medskip


\medskip

\centerline{\bf AMS subject classification: 53A10, 53C12, 35J20}

\end{center}

\

\section{Introduction}

Let $S$ be an oriented embedded (or possibly immersed)
hypersurface in a Riemannian manifold $(M^{n+1},g)$. The shape
operator $A_S$ is the symmetric endomorphism of the tangent bundle
of $S$ associated with the second fundamental form of $S$, $b_S$,
by
\[
b_S (X, Y) = g_S ( A_S \, X, Y), \quad \forall \, X , Y \in TS;
\qquad \mbox{here}\qquad g_S = \left. g \right|_{TS}.
\]
The eigenvalues $\kappa_i$ of the shape operator $A_S$ are the
principal curvatures of the hypersurface $S$. The $k$-curvature of
$S$ is define to be the $k$-th symmetric function of the principal
curvatures of $S$, i.e.
\[
H_k (S) : =  \sum_{i_1 < \ldots < i_k} \kappa_{i_1}\ldots
\kappa_{i_k} .
\]
Hence, when $k=1$, $H_1$ is equal to $n$ times the mean curvature
of $S$. When $M^{n+1}=\R^{n+1}$ is the Euclidean space, $H_2$ is
equal to $\frac{n(n-1)}{2}$ times the scalar curvature of $S$ and
$H_n$ is equal to the Gauss-Kronecker curvature of $S$. In this
paper we are interested in the existence of hypersurfaces in
$M^{n+1}$ whose  $k$-curvature is constant. Hypersurfaces with
constant mean curvature, constant scalar curvature or constant
Gauss-Kronecker curvature in Euclidean space or space forms
constitute an important class of submanifolds. In Riemannian
manifolds very few examples of constant $k$-curvature
hypersurfaces are known, except when $k=1$.

\medskip

R. Ye \cite{Ye-1}, \cite{Ye-2} has proved the existence of a local
foliation by constant mean curvature hypersurfaces which
concentrate at a point (which is required to be a nondegenerate
critical point of the scalar curvature function). We extend the
result and methods of \cite{Ye-1} to handle the case $k=2, \ldots,
n$. No extra curvature hypotheses are required. In particular, we
prove the existence of foliations of a neighborhood of any
nondegenerate critical point of the scalar curvature of $(M^{n+1},
g)$ by constant Gauss-Kronecker or constant scalar curvature
hypersurfaces. As in \cite{Ye-1} the idea is to perturb $\bar
S_\rho (p) $, a geodesic sphere with small radius $\rho >0$
centered at a point $p$. A simple computation will show that $\bar
S_\rho (p) $ is close to being a constant k-curvature hypersurface
as $\rho$ tends to $0$ and in fact
\[
\sigma_k  ( \bar S_\rho  (p) ) = C_n^k \, \rho^{-k} + \calO
(\rho^{-k+2}),
\]
In this paper, we show that it is possible to perturb $\bar S_\rho
(p)$ for every small radius, to a constant $k$-curvature
hypersurface equal to $C_n^k \, \rho^{-k}$ for any $1\le k\le
n-1$, provided $p$ is close to a nondegenerate critical point of
the scalar curvature of $M$. The analysis here is inspired from
the one  performed in \cite{Ye-1}. In fact, independently of the
value of $k$, the linearized $k$-curvature operator about the unit
Euclidean sphere is always a multiple of $\Delta_{S^n} +n$, the
linearized mean curvature operator about the unit Euclidean
sphere. This implies that, as in \cite{Ye-1}, to perform the
perturbation of a small geodesic sphere, one has to overcome the
problem of the existence of $(n+1)$-dimensional kernel of
$\Delta_{S^n} +n$, kernel which is related to the invariance of
$k$-curvature with respect to the action of isometries (in the
case of the unit sphere, this kernel is only generated by
translations). This is where, as in \cite{Ye-1} we use the fact
that we are close to a nondegenerate critical point of the scalar
curvature of the ambient manifold.

\medskip

We notice that the analysis performed in  \cite{Ye-1} is specific
to treat the case of mean curvature, namely $k=1$ and,
unfortunately, can't be used to treat the general case $k=2,
\ldots, n$. The main technical result of this paper is a precise
expansion of geometric operators (first and second fundamental
forms) for perturbed geodesic sphere (see
Proposition~\ref{pr:2.1}, Proposition~\ref{pr:3-2} and
Proposition~\ref{pr:3-4}). We believe that these expansions are of
independent interest and can be used in many other construction
\cite{Mal-Pac}. Our main result is~:
\begin{theorem}
Suppose that $p_0$ is a nondegenerate critical point of the scalar
curvature ${\mathcal R}$ of $M$. Then there exists $\rho_0 >0$,
such that for all $\rho\in (0, \rho_0)$, the geodesic sphere $\bar
S_{\rho} (p_0)$ may be perturbed to a constant $k$-curvature
hypersurface $S_{\rho}$ with $H_k = C_n^k \, \rho^{-k}$. Moreover
these $k$-curvature hypersurfaces constitute a local foliation of
a neighborhood of $p_0$. \label{th:existence}
\end{theorem}
The existence of the hypersurfaces is not so difficult and can be
obtained rather easily. The fact that they constitute a local
foliation requires more work. The leaves $S_\rho$ are small
perturbation of geodesic spheres in the sense that $S_\rho$ is a
normal graph over $\bar S_\rho(p_0)$ for some function $\bar
w_\rho$ which is bounded by a constant times $\rho^2$.

\medskip

The hypersurface $S_\rho$ is a small perturbation of $\bar S_\rho
(p_0)$ in the sense that it is the normal graph of some function
(with $L^\infty$ norm bounded by a constant times $\rho^2$) over a
geodesic sphere obtained centered at a point at distance bounded
by a constant times $\rho^2$ of $p_0$.

\

Existence of families of constant mean curvature hypersurfaces
concentrating along positive dimensional limit sets is obtained by
R. Mazzeo and F. Pacard in \cite{Maz-Pac} and then in
collaboration with the author in  \cite{Mah-Maz-Pac} in a more
general setting.

\

 The paper is organized in the following way: In section 2 we
expand the coefficients of the metric in normal geodesic
coordinates. Section 3 will be devoted to the expansion of the
first fondamental form, second fondamental form and the Shape
operator of the perturbed geodesic spheres. Using these, we derive
in section 4, the expansion of the $k$-curvature of the perturbed
spheres. Section 5 is devoted to the proof of the main result of
this paper, theorem \ref{th:existence}.

\section{Expansion of the metric in geodesic normal coordinates}

In this Section we introduce geodesic normal coordinates in a
neighborhood of a point $p \in M$. We choose an orthonormal basis
$E_i$, $i=1, \ldots, n+1$, of $T_pM$.

\medskip

Consider, in a neighborhood of $p$ in $M$, normal geodesic
coordinates
\[
F(x) : = \exp^M_p (x_i\, E_i), \qquad x := (x_{1}, \ldots,
x_{n+1}),
\]
where $\exp^M$ is the exponential map on $M$ and summation over
repeated indices is understood. This yields the coordinate vector
fields $X_i : = F_* (\del_{x_i})$. As usual, the Fermi coordinates
above are defined so that the metric coefficients
\[
g_{ij} = g( X_i , X_j)
\]
equal $\delta_{ij}$ at $p$. We now compute higher terms in the
Taylor expansions of the functions $g_{ij}$. The metric
coefficients at $q := F(x)$ are given in terms of geometric data
at $p : = F(0)$ and $|x|: = (x_1^2+ \ldots + x_{n+1}^2)^{1/2}$.

\medskip

\noindent {\bf Notation} The symbol $\calO(|x|^r)$ indicates an
analytic function such that it and its partial derivatives of any
order, with respect to the vector fields $x^j \, X_i$, are bounded
by a constant times $|x|^r$ in some fixed neighborhood of $0$.

\medskip

We now give the well known expansion for the metric in normal
coordinates  \cite{Lee-Par}, \cite{Sch-Yau}, \cite{Will}, but we
will  briefly recall the proof in the Appendix for completeness.
\begin{proposition}
At the point $q = F(x)$, the following expansions hold
\begin{eqnarray}
g_{ij} &=& \delta_{ij} + \frac{1}{3} \, g ( R (E_k ,E_i)\, E_\ell
, E_j) \, x_k \, x_\ell \\
&+& \frac{1}{6} \, g(\nabla _{E_k} R (E_\ell, E_i) E_m, E_j ) \,
x_k \, x_\ell \, x_m  + {\mathcal O} (|x|^4)\nonumber
\label{eq:3-1}
\end{eqnarray}
where all curvature terms are evaluated at $p$. Convention ever
repeated indices is understood. \label{pr:2.1}
\end{proposition}

\section{Geometry of spheres}

In this Section, we derive expansions as $\rho$ tends to $0$ for
the metric, second fundamental form and mean curvature of the
sphere $\bar S_\rho(p)$ and their perturbations.

\medskip

Fix $\rho > 0$. We use a local parametrization $z \rightarrow
\Theta (z)$ of $S^{n} \subset T_p M$. Now define the map
\[
G(z) :=  F \, \big(\rho \,(1 - w(z)) \, \Theta (z) \big),
\]
and denote its image by $S_\rho(p,w)$, so in particular
$S_{\rho}(p, 0) = \bar S_\rho (p)$. Because of the definition of
these hypersurfaces using the exponential map, various vector
fields we shall use may be regarded either as fields along $S_\rho
(p,w)$ or as vectors of $T_p M$. To help allay this confusion, we
write
\[
\Theta  : = \Theta^j \, E_j  \qquad \qquad \Theta_i  :=
\partial_{z^i} \Theta^j \, E_j.
\]
These are all vectors in the tangent space $T_p M$.  On the other
hand, the vectors
\[
\Upsilon : = \Theta^j \, X_j  \qquad \qquad \Upsilon_i : =
\partial_{z^i} \Theta^j \, X_j
\]
lie in the tangent space $T_q M$, where $q = F(z)$. For brevity,
we also write
\[
w_j :=  \del_{z^j} w , \qquad w_{ij} : =  \del_{z^i}\, \del_{z^j}
w.
\]

In terms of all this notation, the tangent space to $S_\rho(w)$ at
any point is spanned by the vectors
\begin{equation}
Z_j =  G_* (\del_{z^j}) = \rho  \, ((1- w)\, \Upsilon_j - w_j \,
\Upsilon ), \qquad j=1, \ldots, n. \label{eq:defz0zj}
\end{equation}

\begin{center}{\bf Notation for error terms}\end{center}
 The formulas for the various
geometric quantities of $S_\rho (p,w)$ are potentially very
complicated, and so it is important to condense notation as much
as possible. Fortunately, we do not need to know the full
structure of all of these quantities. Because it is so
fundamental, we have isolated the notational conventions we shall
use in this separate subsection.

\medskip

Any expression of the form $L^j (w)$ denotes a linear combination
of the functions $w$ together with its derivatives with respect to
the vector fields $\Theta_i$ up to order $j$. The coefficients are
assumed to be smooth functions on $S^n$ which are bounded by a
constant independent of $\rho \in (0,1)$ and $p \in M$, in
${\mathcal C}^\infty$ topology.

\medskip

Similarly, any expression of the form $Q^j (w)$ denotes a
nonlinear operator in the functions $w$ together with its
derivatives with respect to the vector fields $\Theta_i$ up to
order $j$. Again, the coefficients of the Taylor expansion of the
corresponding differential operator are smooth functions on $S^n$
which are bounded by a constant independent of $\rho \in (0,1)$
and $p \in M$ in the ${\mathcal C}^\infty$ topology. In addition
$Q^j$ vanishes quadratically at $w = 0$.

\medskip

Finally, any term of the form $L^j {\ltimes} Q^k$ will denote any
finite sum of the product of a linear operators $L^j$ with
nonlinear operators $Q^k$.

\medskip

We also agree that any term denoted $\calO (\rho^d)$ is a smooth
function on $S^n$ which is bounded by a constant (independent of
$p$) times $\rho^d$ in the ${\mathcal C}^\infty$ topology.

\medskip

{\large \bf The first fundamental form} The next step is the
computation of the coefficients of the first fundamental form of
$S_\rho (p, w)$. We set $q  : = G(z)$ and $p := G (0)$. We obtain
directly from (\ref{eq:3-1}) that
\begin{equation}
\begin{array}{rcl}
g ( X_i, X_j) & = &  \delta_{ij} + \frac{1}{3} \, g( R(\Theta ,
E_i ) \, \Theta , E_j ) \, \rho^2 \, (1 - w)^2 \\[3mm]
& + & \frac{1}{6} \, g( \nabla_{\Theta} R (\Theta , E_i ) \,
 \Theta , E_j ) \, \rho^3 \, (1 - w)^3 \\[3mm]
& + & \calO (\rho^4) + \rho^4 \, L^0 (w) + \rho^4 \, Q^0 (w).
\end{array}
\label{lems}
\end{equation}
where all the curvature terms are evaluated at $p$. Observe that
we have
\[
g (\Upsilon, \Upsilon ) \equiv  1  \qquad \quad g (\Upsilon,
\Upsilon_j ) \equiv 0
\]
Using these expansions it is easy to obtain the expansion of the
first fundamental form of $S_\rho(p,w)$.
\begin{proposition}
We have
\begin{equation}
\begin{array}{rlllll}
&&\ds \rho^{-2} \, (1-w)^{-2} \, g( Z_{i}, Z_{j}) =  g (\Theta_i ,
\Theta_j) + \frac{1}{3}\, g( R(\Theta , \Theta_i)\,
\Theta , \Theta_j) \, \rho^2 \, (1-w)^2 \\[3mm]
& + & \frac{1}{6} \, g( \nabla_\Theta R( \Theta , \Theta_i)\,
\Theta , \Theta_j) \, \rho^3 \, (1-w)^3  + (1-w)^{-2} \, w_i\, w_j+\calO (\rho^4)  \\[3mm]
& + &  \rho^4 L^0 (w) + \rho^4 \, Q^0 (w).
\end{array}
\label{eq:3-4}
\end{equation}
where all curvature terms are evaluated at $p$. \label{pr:3-2}
\end{proposition}

{\bf \large The normal vector field} Our next task is to
understand the dependence on $w$ of the unit normal $N$ to $S_\rho
(w)$. Define the vector field
\[
\tilde N  : =  - \, \Upsilon + A^j \, Z_{j} ,
\]
and choose the coefficients $A^j$ so that  $\tilde N$ is
orthogonal to all of the  $Z_i$. This leads to a linear system for
$A^j$.
\[
\sum_j A^j \, g(Z_j, Z_i) = - \rho \, w_i
\]
Observe that
\[
g(\tilde N, \tilde N) =  1 + \rho \, \sum_j A_j \, w_j
\]
The unit normal vector field $N$ about $S_\rho(p, w)$ is defined
to be
\begin{equation}
N  : =  \frac{\tilde N}{g(\tilde N, \tilde N)^{1/2}}
\label{eq:3-5}
\end{equation}

{\bf \large The second fundamental form} We now compute the second
fundamental form. To simplify the computations below, we
henceforth assume that, at the point $\Theta (z) \in S^n$,
\begin{equation}
g ( \Theta_i, \Theta_j ) =\delta_{ij} \qquad \mbox{and}\qquad
\overline \nabla_{\Theta_i} \Theta_{j} =0 \label{eq:3-55}, \quad
i,j = 1, \ldots, n
\end{equation}
(where $\overline \nabla$ is the connection on $TS^{n-1}$).

\begin{proposition}
The following expansions hold
\begin{equation}
\begin{array}{rllll}
- g( \nabla_{Z_i} N,  Z_{j} ) & = &  \rho \, (1-w) \, \delta_{ij}
+ \rho \, w_{ij} + \frac{2}{3} \, g( R(\Theta , \Theta_i )
\, \Theta, \Theta_j ) \, \rho^3 \, (1-w)^3 \\[3mm]
& + & \frac{5}{12} \, g( \nabla_\Theta R(\Theta , \Theta_i ) \,
\Theta, \Theta_j ) \, \rho^4 \, (1-w)^4
\\[3mm]
& - & \frac{1}{3} \, \left( g( R ( \nabla w , \, \Theta_i ) \,
\Theta , \Theta_j) + g( R ( \Theta , \, \Theta_i ) \, \nabla w ,
\Theta_j) \right) \, \rho^3 \\[3mm]
& + &  {\mathcal O} ( \rho^5) + \rho^4 \, L^1 (w)  + \rho  \, Q^1
(w) + \rho \, L^2 (w) \ltimes Q^1(w)
\end{array}\label{eq:3-6d}
\end{equation}
where as usual, all curvature terms are computed at the point $p$.
\label{pr:3-4}
\end{proposition}
{\bf Proof~:} We will first obtain the expansion of $g
(\nabla_{Z_i} \tilde N, Z_j)$. To this aim, we compute
\[
\begin{array}{rllll}
&-& g (\nabla_{Z_i} \tilde N, Z_j)\\[3mm]
 & = & g (\nabla_{Z_i}
\Upsilon,
Z_j) - \sum_k g ( \nabla_{Z_i} (A^k \, Z_k) , Z_j) \\[3mm]
& = &  \frac{1}{1-w} g(\nabla_{Z_i} ((1-w)  \Upsilon ), Z_j) +
\frac{1}{1-w} \, w_i \, g(\Upsilon , Z_j) - \sum_k g ( \nabla_{Z_i} (A^k \, Z_k) , Z_j) \\[3mm]
& = &  \frac{1}{1-w} g(\nabla_{Z_i} ((1-w) \, \Upsilon ), Z_j) -
\frac{\rho}{1-w} \, w_i \, w_j - \sum_k g ( \nabla_{Z_i} (A^k \, Z_k) , Z_j) \\[3mm]
\end{array}
\]
Now, recall that
\[
\sum_k A^k \, g(Z_k, Z_j) =  - \rho \, w_j
\]
Hence
\[
\sum_k  g(\nabla_{Z_i} (A^k \, Z_k) , Z_j) =  - \rho \, w_{ij} -
\sum_k A^k \, g(Z_k, \nabla_{Z_i} Z_j)
\]
Using the fact that
\[
2 \, g(Z_k, \nabla_{Z_i} Z_j) = Z_i \,  g(Z_k, Z_j) + Z_j \,
g(Z_k, Z_i ) - Z_k \, g(Z_i, Z_j)
\]
we conclude that
\begin{eqnarray*}
\sum_k  g(\nabla_{Z_i} (A^k  Z_k) , Z_j) &=&  - \rho  w_{ij} -
\frac{1}{2} \sum_k A^k \, \bigg( Z_i   g(Z_k, Z_j) + Z_j g(Z_k,
Z_i )\\& -& Z_k  g(Z_i, Z_j)\bigg)
\end{eqnarray*}
To analyze the term $\nabla_{Z_i} ((1-w) \, \Upsilon )$, let us
revert for the moment and regard $w$ as functions of the
coordinates $z$ and also consider $\rho$ as a variable instead of
just a parameter. Thus we consider
\[
\tilde{F}(\rho ,z ) =  F \big(\rho (1-w(z)) \Theta (z)\big).
\]
The coordinate vector fields $Z_j$ are still equal to $\tilde F_*
(\del_{z_j})$, but now we also have $Z_0 : =  (1-w) \, \Upsilon =
\tilde F_* (\del_\rho)$, which is the identity we wish to use
below. Now, we write
\begin{eqnarray*}
&&g( \nabla_{Z_i} ((1-w) \, \Upsilon ) , Z_j) + g( \nabla_{Z_j}
((1-w) \, \Upsilon ) , Z_i) \\
&=& g( \nabla_{Z_i} Z_0 , Z_j) + g(
\nabla_{Z_j} Z_0 , Z_i)\\
&=& Z_0 \, g(Z_i, Z_j)
\end{eqnarray*}
Collecting the above we have obtained to formula
\[
\begin{array}{rllll}
- g (\nabla_{Z_i} \tilde N, Z_j) & = &  \frac{1}{2(1-w)} \, Z_0 \,
g(Z_i, Z_j) - \frac{1}{1-w} \, \rho \, w_i \, w_j + \rho \, w_{ij} \\[3mm]
& + & \frac{1}{2} \, \sum_k A^k \, \left( Z_i \,  g(Z_k, Z_j) +
Z_j \,  g(Z_k, Z_i ) - Z_k \, g(Z_i, Z_j)\right)
\end{array}
\]
We will now expand the first and last term in this expression.

\medskip

If the coordinates $y$ are chosen so that $ g(\Theta_i, \Theta_j)
= \delta_{ij}$ at the point where we will compute the shape form,
we have, using the result of Proposition \ref{pr:3-2},
\[
\begin{array}{rlll}
 \frac{1}{2 (1-w)} \, Z_0 \, g( Z_{i}, Z_{j}) & = &  \rho \,
(1-w) \, \delta_{ij} + \frac{2}{3}\, g( R(\Theta ,
\Theta_i)\, \Theta , \Theta_j) \, \rho^3 \, (1-w)^3 \\[3mm]
& + & \frac{5}{12} \, g( \nabla_\Theta R(\Theta , \Theta_i)\,
\Theta , \Theta_j) \, \rho^4 \, (1-w)^4 + \frac{1}{1-w} \, \rho\, w_i \, w_j \\[3mm]
& + & \calO (\rho^5) + \rho^5 \, L^0 (w) + \rho^5 \, Q^0 (w).
\end{array}
\]
Using the same Proposition together with the fact that the
coordinates $y$ are chosen so that $ \bar \nabla_{\Theta_i}
\Theta_j = 0$ at the point where we will compute the shape form,
we also have
\[
\begin{array}{rlllll}
&&Z_i \,  g(Z_k, Z_j) + Z_j \, g(Z_k, Z_i ) - Z_k \, g(Z_i, Z_j)=
\\
[3mm]&& \frac{2}{3} \, \left( g (R(\Theta_k, \Theta_i) \,
\Theta_j, \Theta) + g (R(\Theta_k, \Theta_i) \,
\Theta_j, \Theta) \right) \, \rho^4  \\[3mm]
&& +  {\mathcal O} (\rho^5) + \rho^2 \, L^1 (w) + \rho^2 \, Q^1
(w) + \rho^2 \, L^2 (w) \ltimes L^1 (w)
\end{array}
\]
If the coordinates $y$ are chosen so that $ g(\Theta_i, \Theta_j)
= \delta_{ij}$ at the point where we will compute the shape form,
we have the expansion
\[
A^k = - \frac{w_k}{\rho (1-w)^2} + \rho \, L^1(w) + \rho \, Q^1
(w)
\]
collecting the above estimates, we conclude that
\[
\begin{array}{rllll}
- g( \nabla_{Z_i} \tilde N,  Z_{j} ) & = &  \rho \, (1 - w) \,
\delta_{ij} + \rho \, w_{ij} + \frac{2}{3} \, g( R(\Theta ,
\Theta_i
) \, \Theta, \Theta_j ) \, \rho^3 \, (1-w)^3 \\[3mm]
& + & \frac{5}{12} \, g( \nabla_\Theta R(\Theta , \Theta_i ) \,
\Theta, \Theta_j ) \, \rho^4 \, (1-w)^4
\\[3mm]
& - & \frac{1}{3} \, \left( g( R ( \Theta_k , \, \Theta_i ) \,
\Theta , \Theta_j) + g( R ( \Theta , \, \Theta_i ) \, \Theta_k ,
\Theta_j) \right) \, \rho^3 \, w_k \\[3mm]
& + &  {\mathcal O} ( \rho^5) + \rho^4 \, L^1 (w)  + \rho  \, Q^1
(w) + \rho \, L^2 (w) \ltimes Q^1(w)
\end{array}\label{eq:3-6d}
\]
It remains to observe that
\[
g(\tilde N, \tilde N)^{-1/2} =  1 + Q^1 (w)
\]
This finishes the proof of the estimate. \hfill $\Box$

\medskip

{\bf \large The shape operator of perturbed surfaces} Collecting
the estimates of the last subsection we obtain the expansion of
the shape operator of the hypersurface $S_\rho (p, w)$. In the
coordinate system defined in the previous sections, we get
\begin{proposition}
Under the previous hypothesis, the shape operator of the
hypersurface $S_\rho(p,w)$ is given by
\begin{equation}
\begin{array}{rlllllll}
\rho \,  A_{ij} (w) & = & (1+w) \, \delta_{ij} + w_{ij} +
\frac{1}{3} \,g( R ( \Theta , \, \Theta_i ) \, \Theta , \Theta_j)
\, \rho^2 \\[3mm]
& - &  \frac{1}{3}  \, \big[ \, g( R ( \Theta , \, \Theta_i ) \,
\Theta , \Theta_j) \, w + g(R(\Theta, \Theta_i), \Theta, \Theta_k)
\, w_{kj}\\[3mm]
&+& \left( g( R ( \Theta_k  , \, \Theta_i ) \, \Theta , \Theta_j)
+ g( R ( \Theta , \, \Theta_i ) \, \Theta_k , \Theta_j)  \right)
\, w_k  \big] \, \rho^2 \\[3mm]
&+& \frac{1}{4} \, g( \nabla_\Theta R ( \Theta , \,
\Theta_i ) \, \Theta , \Theta_j) \, \rho^3+\calO (\rho^4)  +\rho^3  L^2 (w)  \\[3mm]
& + &   Q^1(w) + L^2 (w) \ltimes L^0(w) + L^2 (w) \ltimes Q^1(w) .
\notag \label{eq:mc}
\end{array}
\end{equation}
where all curvature terms are computed at the point $p$.
\end{proposition}

\section{The $k$-curvature of the perturbed sphere}

Given any symmetric matrix $A$, and any $k =0, \ldots, n$, we
define
\[
\sigma_k (A) : =  \sum_{i_1 < \ldots < i_k} \lambda_{i_1}\ldots
\lambda_{i_k} .
\]
where $\lambda_1, \ldots, \lambda_n$ are the eigenvalues of $A$.
The $k$-th Newton transform of $A$ is defined by
\[
T_k (A) : = \sigma_k(A) \, I - \sigma_{k-1}(A) \, A + \cdots +
(-1)^k \, A^k.
\]
with $T_{n}(A) = 0$. Now suppose that $A = A(t)$ depends smoothly
on a parameter $t$, it is proved in \cite{R} that
\begin{equation}
\frac{d\,}{dt} \sigma_k (A)= \mbox{Tr} \left( T_{k-1}(A) \,
\frac{d\,}{dt} A \right) \label{eq:trinp}
\end{equation}
From this computation, it follows at once that, given any $n
\times n$ symmetric matrix $H$,
\[
\sigma_k (I+H) = C_n^k+ C_{n-1}^{k-1} \, \mbox{Tr} (H) + {\mathcal
O} (|H|^2)
\]
Using this together with the previous expansion of the shape
operator, it is not hard to check that the $k$-curvature of the
hypersurface $S_\rho (p, w)$ can be expanded as
\[
\begin{array}{rllll}
&&\rho^{k}  \, H_k (S_\rho (p,w) ) \\[3mm]
& = & C_n^k + C_{n-1}^{k-1}
\, \left[ (\Delta_{S^n}+n) \, w - \frac{1}{3} \, \mbox{Ric}
(\Theta, \Theta) \, \rho^2  - \frac{1}{4} \, \nabla_\Theta
\mbox{Ric}
(\Theta, \Theta ) \, \rho^3 \right. \\[3mm]
& + & \frac{1}{3} \, \left( \mbox{Ric}  (\Theta, \Theta) + 2 \,
\mbox{Ric} (\nabla  \cdot , \Theta ) - g(R(\Theta,
\nabla  \cdot ) \Theta, \nabla \cdot ) \right) \, w \, \rho^2 \\[3mm]
& + & \left. {\mathcal O}(\rho^4) + \rho^3 \, L^2 (w) + Q^1(w) +
L^2 (w) \ltimes L^0(w) + L^2 (w) \ltimes Q^1(w) \right]
\end{array}
\]
where as usual, all curvature terms are computed at $p$. Here we
have defined
\[
\mbox{Ric} (\nabla  \cdot , \Theta ) : = \mbox{Ric} (e_i , \Theta
) \, e_i
\]
and
\[
g(R(\Theta, \nabla  \cdot ), \Theta, \nabla \cdot )  :=
g(R(\Theta, e_i), \Theta, e_j ) \, e_i \, e_j
\]
if $e_1, \ldots, e_n$ is an orthonormal frame field of $T_{\bar q}
S^n$ satisfying $\bar \nabla_{e_i} \, e_j = 0$ at the point $\bar
q \in S^n$ where these expressions are computed. It will be
convenient to set
\[
{\mathcal L} : = \frac{1}{3} \, \left( \mbox{Ric}  (\Theta,
\Theta) + 2 \, \mbox{Ric} (\nabla , \Theta ) - g(R(\Theta, \nabla
), \Theta, \nabla  ) \right)
\]

Now observe that a similar expansion is valid in Euclidean space
and in this case the expansion of $\rho^{-k} \, H_k (p,w)$ does
not depend on $\rho$ (nor on $p$). This means that the nonlinear
operator
\[
Q^2 : = Q^1+ L^2 \ltimes L^0 + L^2 \ltimes Q^1
\]
can be decomposed into its value in Euclidean space and a similar
operator all of whose coefficients are bounded by $\rho$. This
fact can also be recovered by going through all the above
expansions. Therefore, we can write
\[
Q^2 = Q_e^2 + \rho \, Q_r^2
\]
where $Q_e^2$ is the corresponding nonlinear operator when the
metric is Euclidean and hence it does not depend on $\rho$; while
$ \rho \, Q_r^2$ denotes the discrepancy induced by the curvature
of the metric $g$ on $M$. Both $Q_e^2$ and $Q_r^2$ satisfy the
usual properties.

\section{Existence of foliations by constant $k$-curvature
hypersurfaces}

Assume that we are given $p_0 \in M$, a nondegenerate critical
point of the scalar curvature ${\mathcal R}$ on $M$. We would like
to find a small function $w \in {\mathcal C}^{2, \alpha} (S^n)$
and a point $p$ close to $p_0$ such that
\[
H_k (S_\rho (p,w)) = C_n^k \, \rho^{-k}
\]
In view of the previous expansion, this amount to solve the
nonlinear equation
\begin{equation}
\begin{aligned}
(\Delta_{S^n}+n) \, w &= \frac{1}{3} \, \mbox{Ric} (\Theta,
\Theta) \, \rho^2  + \frac{1}{4} \, \nabla_\Theta \mbox{Ric}
(\Theta, \Theta ) \, \rho^3 - {\mathcal O}(\rho^4) \\[3mm]
& -  \rho^2 \, {\mathcal L} \, w - \rho^3 \, L^2 (w) - Q^2(w)
\end{aligned} \label{eq:5.1}
\end{equation}

We denote by $\Pi$ and $\Pi^\perp$ the $L^2$-orthogonal
projections of $L^2(S^n)$ onto $\mbox{Ker} (\Delta_{S^n} + n)$ and
$\mbox{Ker} (\Delta_{S^n} + n)^\perp$, respectively. Recall that
the kernel of $\Delta_{S^n} + n$ is spanned by $\varphi_i$, for
$i=1, \ldots, n+1$, the restriction to the unit sphere of $x_i$,
the coordinates functions in ${\mathbb R}^{n+1}$.

\medskip

{\bf First fixed point argument} From now on, we assume that the
function $w \in {\mathcal C}^{2, \alpha} (S^n)$ is
$L^2$-orthogonal to $\mbox{Ker} (\Delta_{S^n} + n)$ and we project
the equation (\ref{eq:5.1}) over $\mbox{Ker} (\Delta_{S^n} +
n)^\perp$. We obtain
\begin{eqnarray*}
(\Delta_{S^n}+n) \, w = &\Pi^\perp \, \left[ \frac{1}{3} \,
\mbox{Ric} (\Theta, \Theta) \, \rho^2  + \frac{1}{4} \,
\nabla_\Theta \mbox{Ric} (\Theta, \Theta ) \, \rho^3  - {\mathcal
O}(\rho^4) \right. \\[3mm]
&\left. - \rho^2 \, {\mathcal L} \, w - \rho^3 \, L^2 (w) -Q^2(w)
\right]
\end{eqnarray*}

We define $w_0 \in \mbox{Ker} (\Delta_{S^n} + n)^\perp$ to be the
unique solution of
\begin{equation}
(\Delta_{S^n}+n)\, w_0 = \frac{1}{3} \, \mbox{Ric} (\Theta,
\Theta) \label{eq:5.10}
\end{equation} since $\Pi^\perp \, ( \mbox{Ric}
(\Theta, \Theta)) = \mbox{Ric} (\Theta, \Theta)$. Similarly, we
define $w_1 \in \mbox{Ker} (\Delta_{S^n} + n)^\perp$ to be the
unique solution of
\[
(\Delta_{S^n}+n)\, w_1 =  \frac{1}{4} \, \Pi^{\perp} \left[
\nabla_\Theta \mbox{Ric} (\Theta, \Theta )\right]
\]

It is easy to rephrase the solvability of the nonlinear equation
(\ref{eq:5.1}) as a fixed point problem since the operator
$\Delta_{S^n}+n$ is invertible from the space of ${\mathcal C}^{2,
\alpha}(S^n)$ functions which are $L^2$-orthogonal to $\mbox{Ker}
(\Delta_{S^n}+n)$ into the space of ${\mathcal C}^{0,
\alpha}(S^n)$ functions which are $L^2$-orthogonal to $\mbox{Ker}
(\Delta_{S^n}+n)$. We write $w : = \rho^2 \, w_0  + \rho^3 \, w_1+
 \rho^4 \, v$, so that it remains to solve an equation which can be written for
short as
\[
(\Delta_{S^n}+n) \, v = - {\mathcal O} (1) - \rho^{-2} \,
{\mathcal L} w - \rho^{-1} \, L^2 (w) - \rho^{-4} \, Q^2 (w)
\]
Applying a standard fixed point theorem for contraction mappings,
it is easy to check that there exists a constant $\kappa > 0$,
which is independent of the choice of the point $p \in M$, such
that there exists a unique fixed point in ball of radius $\kappa$
in ${\mathcal C}^{2, \alpha} (S^n)$, provided $\rho$ is chosen
small enough, say $\rho \in (0, \rho_0)$. We denote by $v_p$ this
solution and define
\[
w_p : = \rho^2 \, w_0 + \rho^3 \, w_1 + \rho^4 \, v_p.
\]
It is easy to check that, reducing the value of $\rho_0$ if this
is necessary,
\begin{equation}
\| w_p - w_{p'} \|_{{\mathcal C}^{2, \alpha}(S^n)} \leq c\, \rho^2
\, \mbox{dist} (p,p'), \label{eq:5.2}
\end{equation}
for some constant $c$ which does not depend on $\rho \in (0,
\rho_0)$ nor on $p$ or $p'$. In addition, the mapping
\[
(\rho, p) \in (0, \rho_0) \times M \longrightarrow w_p \in
{\mathcal C}^{2, \alpha}(S^n)
\]
is smooth and
\[
\| D_p w_p \|_{{\mathcal C}^{2, \alpha}(S^n)} + \rho \, \|
\partial_\rho w_p \|_{{\mathcal C}^{2, \alpha}(S^n)} \leq c \,
\rho^2
\]
for some constant $c$ which does not depend on $\rho \in (0,
\rho_0)$ nor on $p$.

\medskip

{\bf Second fixed point argument} It now remains to project the
equation (\ref{eq:5.1}) where $w$ has been replaced by $w_p$, over
$\mbox{Ker} (\Delta_{S^n} +n)$. To this aim, we recall the nice
and key observation from \cite{Ye-1}.

\medskip

\noindent  The  problem is to compute the $L^2$-projection of the
term $g( \nabla_\Theta R ( \Theta ,  \Theta_i )  \Theta ,
\Theta_j)$ over the kernel of the operator $\Delta_{S^n}+ n$. This
amounts to compute, for any $m=1, \ldots, n+1$, the quantity
\[
B_m : = \sum_{i,j,k,\ell} \, g(  \nabla_{E_j} R (E_i,
E_k) E_i, E_\ell) \, \int_{S^n} x_j \, x_k\, x_\ell \, x_m \\[3mm]
\]
Now to evaluate this quantity, simply use the fact that the
integral vanishes unless all indices are all equal or constitute
two pairs of equal indices. Using this, together with the
symmetries of the curvature tensor which imply that $R(E,E)=0$, we
obtain
\[
\begin{array}{rllll}
B_m & = & g(\nabla_{E_m} R (E_i, E_m) E_i, E_m)  \,
\left(\int_{S^n}
x_1^4 - 3 \, \int_{S^n} x_1^2 \, x_2^2\right)\\[3mm]
& + & g(\nabla_{E_m} R (E_i, E_j) E_i + 2 \nabla_{E_j} R (E_i,
E_m) E_i, E_j)  \, \int_{S^n} x_1^2 \, x_2^2
\end{array}
\]
Now, use second Bianchi identity
\[
g(\nabla_{E_m} R (E_i, E_j) E_i, E_j) = 2 \, g(\nabla_{E_j} R
(E_i, E_m) E_i, E_j)
\]
together with the fact that
\[
\int_{S^n} x_1^4 = 3 \, \int_{S^n} x_1^2 \, x_2^2 =
\frac{3}{(n+3)} \, \int_{S^n} x_1^2
\]
To conclude that
\[
\Pi \, (g( \nabla_\Theta R ( \Theta , \, \Theta_i ) \, \Theta ,
\Theta_j)) =  - \frac{1}{n+3} \, g( \nabla {\mathcal R} , x_i \,
E_i)
\]
where ${\mathcal R}$ denotes the scalar curvature function,
computed at $p$.

\medskip

Therefore, the projection of the equation (\ref{eq:5.1}) over
$\mbox{Ker } (\Delta_{S^n}+n)$ yields
\[
 g( \nabla {\mathcal R} , x_i \, E_i) =  V_p
 \]
 where we have defined
 \[
 V_p : = 4 \, (n+3) \, \Pi \,
\left[ \rho^{-3} \, {\mathcal O}(\rho^4) + \rho^{-1} \, {\mathcal
L} \, w_p + L^2 (w_p) +  \rho^{-3} \, Q^2(w_p)\right]
\]
Now, using the fact that $p_0$ is a nondegenerate critical point
of the scalar curvature, we conclude easily (applying for example
a topological degree argument) that there exists $p$ close to
$p_0$ satisfying (\ref{eq:5.2}) provided $\rho$ is close enough to
$0$. This gives the existence of constant $k$-curvature leaves for
all $\rho$ small enough, unfortunately it turns out that the point
$p$ is at most at distance a constant times $\rho$ from $p_0$ and
this is not enough to show that the constant $k$-curvature leaves
form a foliation of a neighborhood of $p_0$.

\medskip

To improve this estimate, many observations are due. First,
observe that we can decompose ${\mathcal O}(\rho^4)$ into the sum
of two functions, one of which is homogeneous of degree $4$ (in
the coordinate functions $x_i$) and the other one which is bounded
by a constant times $\rho^5$. The $L^2$-projection of the
homogeneous function of degree $4$ is equal to $0$ since this
homogeneous function is invariant under the change of coordinates
$\Theta$ into $-\Theta$. Hence we conclude that
\[
| \Pi \, ({\mathcal O}(\rho^4) ) | \leq c \, \rho^5
\]
Similarly, observe that $w_0$ and hence ${\mathcal L} \, w_0$ are
invariant under the change $\Theta$ into $-\Theta$ and hence the
$L^2$ projection of ${\mathcal L} \, w_0$ over $\mbox{Ker}
(\Delta_{S^n}+n)$ again identically equal to $0$. Therefore, we
conclude that
\[
| \Pi \, ( {\mathcal L} \, w_p )  | \leq c \, \rho^3
\]
Finally, we use the observation at the end of \S 4. Since the
nonlinear operator $Q^2_e$ preserves functions which are invariant
under the action of $-I$, we conclude that $\Pi (Q^2_e (\rho^2 \,
w_0))=0$ and hence
\[
| \Pi (Q^2 (w_p)) | \leq c \, \rho^5
\]
These precise estimates imply that,
\[
|V_p|\leq c \, \rho^2
\]
for some constant which does not depend on $p$ nor on $\rho$. With
slightly more work, we get using similar arguments that
\begin{equation}
| \Pi \, \left( V_p - V_{p'} \right) | \leq c \, \rho^2 \,
\mbox{dist} (p,p') \label{eq:5.12}
\end{equation}
Now, for all $\rho$ small enough, we can find a solution of
(\ref{eq:5.1}) using a fixed point argument for contraction
mapping, in the geodesic ball of radius $2 \, \rho^2$ centered at
any nondegenerate critical point of ${\mathcal R}$. Moreover, the
solution $p_\rho$ depends smoothly on $\rho$ and
\[
|\del_\rho p_\rho |\leq c\, \rho
\]
This later fact, together with (\ref{eq:5.12}) shows that the
solutions constitute a local foliation. This completes the proof
of the main result.  \hfill{$\Box$}

\medskip

Having derived such precise estimates, we can compute the
expansion of the $n$-dimensional volume of the leaves of the
foliation as well as the $(n+1)$-dimensional volume enclosed by
each leaf.

\begin{proposition}
For all $\rho$ small enough the following expansions hold for the
$n$-dimensional volume of $S_\rho$
\[
\mbox{Vol}_n (S_\rho) = \rho^n \, \mbox{Vol}_n (S^n) \, \left( 1-
\frac{1}{2(n+1)} \, {\mathcal R} \, \rho^2 + {\mathcal O} (\rho^4)
\right)
\]
and the $(n+1)$-dimensional volume of the set $B_\rho$ enclosed by
$S_\rho$ and containing the point $p_0$
\[
\mbox{Vol}_{n+1} (B_\rho) = \frac{1}{n+1} \, \rho^{n+1} \,
\mbox{Vol}_n (S^n) \, \left( 1 - \frac{n+2}{2n (n+3)} \, {\mathcal
R} \, \rho^2 + {\mathcal O} (\rho^4)\right)
\]
where the scalar curvature is computed at $p_0$, a nondegenerate
critical point of ${\mathcal R}$.
\end{proposition}
{\bf Proof~:} Integrating (\ref{eq:5.10}) over $S^n$ we find
\[
n \, \int_{S^n} w_0 = \frac{1}{3} \,  \int_{S^n} \mbox{Ric}
(\Theta, \Theta)
\]
Now, plugging the expansion of $w_p$ into the expression of the
first fundamental form given in Proposition~\ref{pr:3-2}, we find
the expansion of $h$ the induced metric on $S_\rho$
\begin{equation}
\begin{aligned} \ds \rho^{-2} \, h_{ij}& = (1 - 2 \rho^2 w_0 - 2 \rho^3 w_1 )
\delta_{ij} + \frac{1}{3}\, g( R(\Theta , \Theta_i)\, \Theta ,
\Theta_j) \, \rho^2 \\[3mm]
&+ \frac{1}{6} \, g( \nabla_\Theta R( \Theta , \Theta_i)\, \Theta
, \Theta_j) \, \rho^3 +  \calO (\rho^4)
\end{aligned}
\end{equation}
This implies that
\[
\rho^{-n} \sqrt |h| = 1 - n \, \rho^2 \, w_0 -n \, \rho^3 \, w_1
-\frac{1}{6} \, \mbox{Ric} (\Theta, \Theta) \, \rho^2 -
\frac{1}{12} \, \nabla_\Theta \mbox{Ric}( \Theta , \Theta) \,
\rho^3 + {\mathcal O} (\rho^4)
\]
The first estimate follows from integrating this expansion using
the fact that the integral of $w_1$ and the integral
$\nabla_\Theta \mbox{Ric} \, ( \Theta , \Theta)$ over $S^n$ vanish
together with the fact that~
\[ \int_{S^n} \mbox{Ric}(\Theta,
\Theta)=\frac{1}{n+1}\,\mbox{Vol}_n(S^n)\,{\mathcal R}
\]

\medskip

Next, we consider polar geodesic normal coordinates $(r, \Theta)$
centered at $p_\rho$. In these coordinates the metric $g$ expanded
as~
\begin{equation}
\ds r^{-2} \, g_{ij}  =   \delta_{ij} + \frac{1}{3}\, g( R(\Theta
, \Theta_i)\, \Theta , \Theta_j) \, r^2 + \frac{1}{6} \, g(
\nabla_\Theta R( \Theta , \Theta_i)\, \Theta , \Theta_j) \, r^3 +
\calO (r^4).
\end{equation}
then, the volume form can be expanded as
\[
r^{-n} \, \sqrt |g| = 1 - \frac{1}{6} \, \mbox{Ric} (\Theta,
\Theta) \, r^2 - \frac{1}{12} \, \nabla_\Theta \mbox{Ric} (\Theta,
\Theta) \, r^3+ {\mathcal O} (r^4)
\]
Integration over the set $r \leq \rho \, (1-w_p)$ give
\[
\begin{aligned}
&\mbox{Vol}_{n+1}(B_\rho)\\[3mm]
&=\int\int_{r\le
\rho\,(1-w_p)}\,r^n\,\left(1 - \frac{1}{6} \, \mbox{Ric} (\Theta,
\Theta) \, r^2 - \frac{1}{12} \, \nabla_\Theta \mbox{Ric} (\Theta,
\Theta) \, r^3+ {\mathcal O}
(r^4)\right)\\[3mm]
&=\frac{1}{n+1}\,\rho^{n+1}\,\int_{S^n}(\,1-(n+1)\rho^2\,w_0\,)+{\mathcal O}(\rho^{n+5})\\[3mm]
&-\frac{1}{6}\,\frac{1}{n+3}\,\rho^{n+3}\,
\int_{S^n}(\,1-(n+3)\rho^2\,w_0\,)\,\mbox{Ric}( \Theta
,\Theta)\\[3mm]
&=\frac{1}{n+1}\rho^{n+1}\mbox{Vol}_n(S^n)-\rho^{n+3}\int_{S^n}\,w_0-\frac{1}{6}\,\frac{\rho^{n+3}}{n+3}
\int_{S^n}\mbox{Ric}(
\Theta , \Theta)+{\mathcal O}(\rho^{n+5})\\[3mm]
&=\frac{1}{n+1}\rho^{n+1}\mbox{Vol}_n(S^n)-(\,\frac{1}{3n}+\frac{1}{6}\,\frac{1}{n+3})\rho^{n+3}\int_{S^n}\mbox{Ric}(
\Theta , \Theta)+{\mathcal O}(\rho^{n+5})\\[3mm]
&=\frac{1}{n+1}\,\rho^{n+1}\,\mbox{Vol}_n(S^n)\,\left(1 -
\frac{n+2}{2n (n+3)} \, {\mathcal R} \, \rho^2 + {\mathcal O}
(\rho^4)\right)
\end{aligned}
\]
This gives the second estimate. This proves the desired result.
\hfill $\Box$

\section{Appendix : proof of Proposition~\ref{pr:2.1}}

The aim of this Section is to prove Proposition \ref{pr:2.1}.
Observe first that the curve $s \rightarrow \exp^M_p (sE)$ is a
geodesic. Therefore, if $X$ is the unit tangent vector to the
curve we have $\nabla_X X = 0$. Hence we also have $(\nabla_X)^m
X=0$ for all $m \geq 1$. In particular, we have, at $p$,
\[
(\nabla_E)^m  E = 0
\]
for all $E \in T_pM$ and for all $m \geq 1$.

\medskip

Observe that $X_a$ are coordinate vector fields hence
\[
\nabla_{X_a} X_b = \nabla_{X_b} X_a
\]

Taking $E = E_a + \e \, E_b$ and looking for the coefficient of
$\e$ in $\nabla_E E = 0$, we get
\[
\nabla_{E_a} E_b = 0
\]
Looking at the coefficient of $\e$ in $\nabla_E^2 E=0$, we get
\begin{equation}
 2 \, \nabla_{E_a}^2 E_b + \nabla_{E_b} \nabla_{E_a} E_a =0
\label{eq:6.13} \end{equation} Finally, looking at the coefficient
of $\e$ in $\nabla_E^3 E=0$, we get~
 \begin{equation} 2 \, \nabla_{E_a}^3
E_b + (\nabla_{E_b} \nabla_{E_a} + \nabla_{E_a} \nabla_{E_b}  )
\nabla_{E_a} E_a =0 \label{eq:6.14}
\end{equation}
Recall that, by definition
\[
\nabla_X \, \nabla_Y : = R(X,Y) + \nabla_Y \, \nabla_X  +
\nabla_{[X,Y]}
\]
Hence, if $X$ and $Y$ are coordinate vector fields we simply have
\begin{equation}
\nabla_X \, \nabla_Y  X : = R(X,Y)X  + \nabla_Y \, \nabla_X X
\label{eq:6.15}
\end{equation}
 We also have
\begin{equation}
\begin{array}{rllll}
\nabla_Y \, \nabla_X \, \nabla_Y  X & : = & \nabla_Y R (X,Y) X  +
R (\nabla_Y X,Y) X + R (X,\nabla_Y Y) X
\\[3mm]
& + & R (X,Y)  \nabla_Y X+ \nabla_Y^2 \nabla_X X + \nabla_Y
\nabla_{[X,Y]} X
\end{array}
\label{eq:6.16}
\end{equation}
Now use (\ref{eq:6.13}) and (\ref{eq:6.15}) to obtain
\begin{equation}
3 \, \nabla_{E_a}^2 E_b  =  R(E_a,E_b) \, E_a, \label{eq:6.17}
\end{equation}
Similarly, use (\ref{eq:6.14}) and (\ref{eq:6.16}) to obtain
\begin{equation}
2 \, \nabla_{E_a}^3 \, E_b + R (E_b,E_a) \, \nabla_{E_a} E_a + 2
\nabla_{E_a} \nabla_{E_b}  \nabla_{E_a} E_a =0\label{eq:6.17}
\end{equation}
Since $\nabla_{E_a} E_b=0$, we get
\[
2 \, \nabla_{E_a}^3 \, E_b +  2 \nabla_{E_a} \nabla_{E_b}
\nabla_{E_a} E_a=0
\]
Using this, we conclude that
\[
\begin{aligned}
2 \, \nabla_{E_a}^3 \, E_b = -2\,\nabla_{E_a}
\nabla_{E_b}\nabla_{E_a} E_a&=-2\,\nabla_{E_a}(R ( E_b ,
E_a) \, E_a+\nabla_{E_a}\nabla_{E_a} E_b)\\[3mm]
&=-2\,\nabla_{E_a}(R ( E_b , E_a) \, E_a)-2 \, \nabla_{E_a}^3 \,
E_b
\end{aligned}
\]
Hence~
\begin{equation}
2 \, \nabla_{E_a}^3 \, E_b =- \nabla_{E_a} R ( E_b , E_a) \, E_a
\label{eq:6.19}
\end{equation}
 Now, we have~
\[
X_c \, g_{ab} = g(\nabla_{X_c} X_a , X_b) + g(X_a, \nabla_{X_c}
X_b),
\]
and we get $\left. X_c \, g_{ab} \right|_p =0$. This yields the
first order Taylor expansion
\[
g_{ab} = \delta_{ab} + {\mathcal O} (|x|^2) ,
\]
To compute the second order terms, it suffices to compute $X_c^2\,
g_{ab}$ at $p$ and polarize. We compute
\[
X_c^2 \, g_{ab} = g ( \nabla_{X_c}^2 X_a , X_b) + g( X_a,
\nabla_{X_c}^2 X_b) + 2 \, g( \nabla_{X_c} X_a, \nabla_{X_c} X_b)
\]
Using (\ref{eq:6.17}) we get
\[
\left. X_c^2 \, g_{ab} \right|_p = \frac{2}{3} \, g( R(E_c, E_a)
\, E_c , E_b).
\]
The formula for the second order Taylor coefficient for $g_{ab}$
now follows at once.

\medskip

Similarly, we compute
\begin{eqnarray*}
\left. X^3_c \, g_{ab} \right|_p &=&  g ( \nabla_{X_c}^3 X_a ,
X_b)
+ 3 g( \nabla_{X_c}^2 X_a, \nabla_{X_c} X_b) \\[3mm]
&+& 3 g( \nabla_{X_c} X_a, \nabla_{X_c}^2 X_b) + g( X_a,
\nabla_{X_c}^3 X_b)
\end{eqnarray*}
and using (\ref{eq:6.19}) this gives
\[
\left. X_c^3 \, g_{ab} \right|_p =  g( \nabla_{E_c} R (E_a, E_c)
\, E_b , E_c).
\]
the formula for the second order Taylor expansion for $g_{ab}$
holds at once. \hfill $\Box$

\end{document}